\documentstyle{article}
\begin{document}

\title {\large\bf Solution to the Generalized Champagne Problem on
Simultaneous Stabilization of Linear Systems
\thanks{Supported by National Natural Science Foundation of China (Nos. 60572056, 60528007, 60334020, 60204006,
10471044, and 10372002), National Key Basic Research and Development
Program (Nos. 2005CB321902, 2004CB318003, 2002CB312200), the
Overseas Outstanding Young Researcher Foundation of Chinese Academy
of Sciences and the Program of National Key Laboratory of
Intelligent Technology and Systems of Tsinghua University.}}

\vskip 0.8cm

\author{\normalsize Qiang Guan$^{1)}$\quad Long Wang$^{2)}$ \quad Bican Xia$^{3)}$\quad Lu Yang$^{4)}$\quad Wensheng Yu$^{1),5)}$\quad Zhenbing Zeng$^{4)}$  \\[0.5cm]
\small\it  1)Laboratory of Complex Systems and Intelligence Science, Institute of Automation,\\
\small\it  Chinese Academy of Sciences, Beijing 100080, P. R. China.\\[0.05cm]
\small\it  2)Center for Systems and Control, Department of Mechanics and Engineering Science,\\
\small\it  Peking University, Beijing 100871, P. R. China.\\[0.05cm]
\small\it  3)LMAM and School of Mathematical Sciences, Peking University\\
\small\it   Beijing 100871, P. R. China.\\[0.05cm]
\small\it  4)Shanghai Institute of Theoretical Computing, Software Engineering Institute,\\
\small\it  East China Normal University, Shanghai 200062, P. R. China.\\[0.05cm]
\small\it  5)National Key Laboratory of Intelligent Technology and Systems, \\
\small\it  Tsinghua University, Beijing 100084, P. R. China.\\[0.1cm]
}
\date{} \maketitle

\begin{quote} \begin{small} \noindent \footnotesize
{\bf Abstract:} {\it The well-known ``Generalized Champagne
Problem" on simultaneous stabilization of linear systems is solved
by using complex analysis \cite{A73,C78,G69,N52,R87} and Blondel's
technique \cite{B94,BG93,BGMR94}. We give a complete answer to the
open problem proposed by Patel et al. \cite{P99,PDV02}, which
automatically includes the solution to the original ``Champagne
Problem" \cite{BG93,BGMR94,BSVW99,LKZ99,P99,PDV02}. Based on the
recent development in automated inequality-type theorem proving
\cite{Y98,Y99,YHX01,YX05,YX03}, a new stabilizing controller
design method is established. Our numerical examples significantly
improve the relevant results in the literature
\cite{LKZ99,PDV02}.}
\\ [0.2cm]
{\bf Keywords:}
 {\it linear systems, stabilization, simultaneous stabilization, champagne
 problem, generalized champagne problem, complex analysis, inequality-type theorem, automated theorem proving. }
 \\ \end{small}
\end{quote}

\vskip 0.5cm \section{ Introduction} \vskip 0.5cm

Simultaneous stabilization of linear systems is a fundamental
issue in system and control theory, and is of  theoretical as well
as practical significance \cite{B94,DFT91,V85}. The basic
statement of the simultaneous stabilization problem of linear
systems is as follows \cite{B94,DFT91,SM82,V85}:

{\it Let $p_1, p_2, \cdots,p_k$ be $k$ scalar linear
time-invariant systems. Under what condition does there exist a
fixed controller $c$ that is stabilizing for each
$p_i(i=1,\cdots,k)$?}

If $k=1$, the above problem is reduced to the stabilization of a
single system. There always exists a stabilizing controller for a
single system provided no unstable pole-zero cancellations occur
\cite{B94,DFT91,V85}. Moreover, once a stabilizing controller of a
single system is found, it is easy to parametrize the set of all
stabilizing controllers of this system
\cite{B94,DLMS80,DFT91,V85}. This parametrization is known as
Youla-Kucera parametrization discovered by Youla et al.
 \cite{YBJ76,YJB76} and Kucera \cite{K79} respectively.

When $k=2$, by using the Youla-Kucera parametrization, it is
possible to rephrase simultaneous stabilization of two systems
into strong stabilization (stabilization with a stable controller)
of a single system. This relationship was discovered for scaler
systems by Saeks and Murray \cite{SM82}, and for multi-variable
systems \cite{B94,DFT91,V85} by Vidyasagar and Viswanadham
\cite{VV82}. For the strong stabilization problem of a single
system, Youla et al. \cite{YBL74} established an elegant
criterion: a system is stabilizable by a stable controller if and
only if it has an even number of real unstable poles between each
pair of real unstable zeros! This remarkable and easily testable
condition is called ``parity interlacing property"
\cite{B94,DFT91,V85,YBL74}.

For $k\geq3$, simultaneous stabilization problem is much more
complicated than one expected \cite{B94,BG93,BGMR92,BGMR94}.
Vidyasagar and Viswanadham \cite{VV82} stated that it is possible
to transform simultaneous stabilization of $k$ systems to strong
stabilization of corresponding $k-1$ systems \cite{B94,DFT91,V85}.
Blondel et al. \cite{BGMR94} proved that simultaneous
stabilization of $k$ systems is equivalent to bistable
stabilization of associated $k-2$ systems. Bistable stabilization
means stabilization with a stable and inverse-stable controller.
Such a controller is called bistable controller or unit controller
\cite{B94}. Although many necessary or sufficient conditions for
simultaneous stabilization of three or more systems were obtained
in recent years (see \cite{B94} and references therein), easily
testable necessary and sufficient conditions have not been found
yet. Blondel and Gevers \cite{BG93} showed that the simultaneous
stabilization of three systems is not rationally decidable
\cite{B94,BG93}, i.e., it is not possible to find tractable
necessary and sufficient conditions for simultaneous stabilization
of three systems that involve only a combination of finite
arithmetical operations (addition, substraction, multiplication
and division), logical operations (`and' and `or') and sign test
operations (equal to, greater than, greater than or equal to, less
than, less than or equal to) on the coefficients of the three
systems!

To illustrate the complexity of the simultaneous stabilization
problem of three systems, Blondel et al. \cite{BG93,BGMR94,BSVW99}
proposed a specific simultaneous stabilization problem called
``Champagne Problem" where the three systems are explicitly given.
Patel \cite{P99} solved this problem by showing that there does not
exist a stabilizing controller. Furthermore, a more general
simultaneous stabilization problem, the ``Generalized Champagne
Problem", was proposed by Patel et al. \cite{P99,PDV02}, which
contains the ``Champagne Problem" as a special case.

The ``Generalized Champagne Problem " is essentially concerned
with determining the exact range of a specific parameter in
simultaneous stabilization problem \cite{P99,PDV02}.  Up to now,
there are only some numerical estimates on this range in the
literature \cite{LKZ99,PDV02}.

It should be pointed out that the ``Champagne Problem" and the
``Generalized Champagne Problem" are not of much engineering
significance in themselves since they both deal with some specific
systems. However, by studying these specific problems, we can
develop new tools and new methods for general simultaneous
stabilization problem. As a matter of fact, by studying some
specific systems, Blondel et al. revealed the inherent relationship
\cite{B94,BG93,BGMR94,P99} between the simultaneous stabilization
problem and complex analysis theory \cite{A73,C78,G69,N52,R87}, and
derived the deep theoretic results of rational undecidability for
three systems \cite{B94,BG93}.

In this paper, the well-known ``Generalized Champagne Problem" on
simultaneous stabilization of linear systems is solved by using
complex analysis theory \cite{A73,C78,G69,N52,R87} and Blondel's
technique \cite{B94,BG93,BGMR94}. We give a complete answer to the
open problem proposed by Patel et al. \cite{P99,PDV02}, which
automatically includes the solution to the original ``Champagne
Problem". Based on the recent advances in automated
inequality-type theorem proving \cite{Y98,Y99,YHX01,YX03,YX05}, a
novel stabilizing controller design method is established. Our
numerical examples significantly improve the relevant results in
the literature \cite{LKZ99,PDV02}.

The paper is organized as follows. We state the generalized
champagne problem and give a complete theoretical solution in
Section 2. Section 3 deals with the controller design, and some
numerical examples are provided. Finally, the conclusion is
contained in Section 4.

\vskip 0.7cm \section{Solution to the Generalized Champagne
Problem}

In this paper, all polynomials are of real coefficients. We denote
by $P$ the set of polynomials,  $P^n$ the set of $n$-th order
polynomials where $n$ is an non-negative integer, $H$ the set of
Hurwitz stable polynomials (all roots lie within left half of the
complex plane), $H^n$ the set of $n$-th order Hurwitz polynomials,
$S$ the set of Schur stable polynomials (all roots lie outside the
unit circle)\footnote{ Note that, some researchers use the concept
of Schur polynomials for polynomials with all roots lying inside
unit circle. Most papers on simultaneous stabilization adopt the
definition used in this paper, which makes it easier to use the
results in complex analysis. It is only a usage and doesn't affect
the problem in essence.} and $S^n$ the set of $n$-th order Schur
polynomials. The variable of polynomials is $s$ or $z$ (we usually
use $s$ for Hurwitz polynomials and $z$ for Schur polynomials). Let
$C$ denote the complex plane, $D=\{z\in C: |z|<1\}$ be the open unit
disc and $\mbox{cl}(D)=\{z\in C: |z|\leq 1\}$ be the closure of $D$,
i.e., the closed unit disc. As usual, $\overline{z}$ stands for the
complex conjugate of $z\in C$, $C_\infty=C\cup \{\infty\}$ the
extended complex plane and $C_{+\infty}=\{s\in C: \mbox{Re}(s)\ge
0\}\cup \{\infty \}$ the extended closed right half complex plane.

Throughout this paper, plants and controllers are restricted to
single-input single-output systems that are described by linear,
time-invariant, real rational transfer functions\footnote{Here, the
real rational transfer functions are not necessarily proper, that is
to say, the degree of numerator may be greater than that of
denominator. In the literature on simultaneous stabilization
\cite{B94,BG93,BGMR94,BSVW99}, controllers do not have to be proper.
It is mainly because that simultaneous stabilization problems can be
transformed into the study on the properties of some analytic
functions on the extend complex plane \cite{A73,C78,G69,N52,R87}. In
addition, Blondel \cite{B94} gave the conclusion that if $k$ plants
are simultaneously stabilizable by a non-proper controller, then
they are also simultaneously stabilizable by a proper controller. In
fact, once a non-proper controller simultaneously stabilizes $k$
systems, since the roots of a polynomial continuously depend on its
coefficients \cite{G59,YZH96}, a proper stabilizing controller can
be obtained with a sufficiently small perturbation imposed on the
denominator polynomial of the original controller. This technique
will also be used in the following section on controller design.}.

Given a system $p$, a controller $c$ stabilizing $p$ means that the
four transfer functions $pc(1+pc)^{-1},c(1+pc)^{-1},p(1+pc)^{-1}$
and $(1+pc)^{-1}$ are all stable, i.e., the denominator polynomials
of the transfer functions are stable (which means Hurwitz stable in
continuous-time case or Schur stable in discrete-time case). The
degree of a controller is the maximum of the degrees of its
denominator and numerator.

We denote by $U(\mbox{cl}(D))=\{c(z)=\displaystyle\frac{y(z)}{x(z)}:
x(z)\in S,y(z) \in S\}$ the set of all bistable or unit controllers
in discrete time case. Obviously, any element in $U(\mbox{cl}(D))$
and its inverse are analytic on $\mbox{cl}(D)$.

The following problem is the well-known ``Champagne Problem"
\cite{BG93,BGMR94, BSVW99,LKZ99,P99,PDV02} of simultaneous
stabilization. It illustrates the difficulty in the simultaneous
stabilization of three systems in general. It also shows that,
simple problems do not always have simple answers.

{\bf Champagne Problem} \cite{BG93,BGMR94, BSVW99,LKZ99,P99,PDV02}
\quad {\it Does there exist a controller that simultaneous
stabilizes the following three plants:
$$p_1(s)=\displaystyle\frac{2}{17}\displaystyle\frac{s-1}{s+1},
\quad p_2(s)=\displaystyle\frac{(s-1)^2}{(9s-8)(s+1)},\quad p_3(s)=0?$$}

Patel \cite{P99} solved this problem in 1999 by showing that there
does not exist a stabilizing controller. More generally, he
considered the following problem:

{\bf Generalized Champagne Problem} \cite{P99,PDV02} \quad {\it
What is the range of $\delta$ for the existence of a controller
that simultaneously stabilizes the following three plants:
$$p_1(s)=\displaystyle\frac{2\delta (s-1)}{s+1},
\quad p_2(s)=\displaystyle\frac{2\delta (s-1)^2}{((1+\delta
)s-(1-\delta ))(s+1)},\quad p_3(s)=0 ? $$ where $\delta $ is a
real number. }

Obviously, ``Generalized Champagne Problem" focuses on determining
the range of $\delta$ when there exists a simultaneously stabilizing
controller while ``Champagne Problem" asks whether $\delta
=\displaystyle\frac{1}{17}$ is in this range.

Regarding ``Generalized Champagne Problem", as discussed in Patel
\cite{P99}, there does not exist such a controller when $\delta
<\displaystyle\frac{1}{16}$ \footnote{In fact, Patel \cite{P99}
admitted $\delta > 0$. Otherwise this condition should be
$0<|\delta| <\displaystyle\frac{1}{16}$. Indeed, from the result in
the sequel, when $\delta =0$ or $\delta
<-\displaystyle\frac{1}{16}$, there does exist a simultaneously
stabilizing controller for ``Generalized Champagne Problem".}. This
conclusion certainly answers ``Champagne Problem" and Patel also
stated that it remains an open problem \cite{P99} whether there
exists such a controller when $\delta =\displaystyle\frac {1}{16}$.
Leizarowitz et al. \cite{LKZ99,PDV02} conjectured the minimum value
of $\delta$ for ``Generalized Champagne Problem" is
$\displaystyle\frac{1}{2e}=\displaystyle\frac {1}{5.4366}$ and they
found a controller for this value of $\delta$. But in 2002, Patel et
al. \cite{PDV02} showed that there exists a stabilizing controller
$c=\displaystyle\frac{y(s)}{x(s)}$ where $x(s)\in H^9,y(s)\in P^9$
when $\delta =\displaystyle\frac{1}{6.719367588932806}$. Note that
the degree of the controller therein equals 9. Hence, the conjecture
mentioned in \cite{LKZ99} is invalid. Thus for ``Generalized
Champagne Problem", it is still an open problem whether there exists
a stabilizing controller when $\delta \in
[\displaystyle\frac{1}{16},\displaystyle\frac{1}{6.719367588932806})$.

{\it  For ``Generalized Champagne Problem" with $\delta
=\displaystyle\frac {1}{16}$, since the roots of a polynomial
continuously  depend on its coefficients \cite{G59,YZH96}, it is
easy to know that there does not exist a controller that
simultaneously stabilizes the three plants. Thus we have given an
answer to the problem proposed in Patal \cite{P99}. }

It should be pointed out that the above results presented in Patel
et al. \cite{PDV02} and Leizarowitz et al. \cite{LKZ99} are both
based on numerical analysis. The degrees of the controllers provided
are high and the conclusions are not theoretically complete.

The following theorem is the main result of this paper which gives a
completely theoretical solution to ``Generalized Champagne Problem".

{\bf Theorem 1} \quad {\it The necessary and sufficient condition
for the existence of a controller that simultaneously stabilizes
the following three plants:
\begin{equation}
p_1(s)=\displaystyle\frac{2\delta (s-1)}{s+1},\quad
p_2(s)=\displaystyle\frac{2\delta (s-1)^2}{((1+\delta )s-(1-\delta
))(s+1)},\quad p_3(s)=0
\end{equation}
where $\delta $ is a real number, is $\delta =0$ or $|\delta|
> \displaystyle\frac {1}{16}.$}

To prove the above theorem, the following two well-known results in
complex analysis will be used.

{\bf Lemma 1 ($\displaystyle\frac{1}{16} - $ Theorem)}
\cite{B94,G69,N52,P99} \footnote{There is a misprint in \cite{P99}
for the reference of this lemma where Condition 2 is misprinted as
$f(0)=1.$} { \quad {\it Any analytic function $f$ on D that
satisfies:

1) $f(z)=0,z\in D \Leftrightarrow z=0,$

2) $f'(0)=1,$

\noindent contains an open ball centered at $z=0$ with radius
$\displaystyle\frac{1}{16}$ in its range on D but not always a
larger ball.}

{\bf Lemma 2 (Runge's Theorem)} \cite{C78,R87} \quad {\it If
$\Omega$ is an open subset of $C$ such that $C_\infty\setminus
\Omega$ is connected then for each analytic function $f$ on $\Omega$
there is a sequence of polynomials $\{q_n\}$ such that $q_n$
uniformly converges to $f$ on any compact subset of $ \Omega $.}
\footnote{This lemma is a special case of Runge's Theorem, i.e., a
corollary of the original Runge's Theorem \cite{C78,R87}.}

We are now in the position to give a proof of Theorem 1.

{\bf Proof of Theorem 1} \quad The well known linear fractional
transformation $s=\displaystyle\frac{1+z}{1-z}$ is a one-to-one
mapping of $\mbox{cl}(D)$ onto $C_{+\infty}$ and its inverse is
$z=\displaystyle\frac{s-1}{s+1}$. By this transformation,
``Generalized Champagne Problem" is equivalent to the simultaneous
stabilization problem of the following three plants in discrete time
case:
\begin{equation}
p_{1d}(z)=2\delta z,\quad p_{2d}(z)=\displaystyle\frac{2\delta
z^2}{z+\delta },\quad p_{3d}(z)=0, \label{dp}
\end{equation}
where $\delta$ is a real number.

The proof contains the following three steps.

{\bf Step 1}\quad Firstly, we prove that the three plants defined by
Eq. (\ref{dp}) are not simultaneously stabilizable in discrete time
case when $0<|\delta| <\displaystyle\frac{1}{16}$.

Assumed by contradiction that, for
$0<|\delta|<\displaystyle\frac{1}{16}$, there is a stabilizing
controller $c(z)$ for the three plants
$p_{1d}(z),p_{2d}(z),p_{3d}(z)$. Then, since $c(z)$ stabilizes
$p_{3d}(z)=0$, $c(z)$ must be stable. Moreover, $c(z)$ also
stabilizes $p_{1d}(z)=2\delta z$ and
$p_{2d}(z)=\displaystyle\frac{2\delta z^2}{z+\delta}$ and we have
the following equations:
\begin{eqnarray}
2\delta zc(z)+1&=&u_1(z)\in  U(\mbox{cl}(D))\\
2\delta z^2c(z)+z+\delta&=&u_2(z)\in U(\mbox{cl}(D))
\end{eqnarray}

Set $$f(z) = 2\delta z^2c(z)+z=z(2\delta zc(z)+1).$$ Obviously,
$f(z)$ is analytic on $D$ and $f(0)=0$ as well as $f'(0)=1$. Also,
we have
\begin{eqnarray}
f(z)&=&z u_1(z)\label{1st}\\
f(z)+\delta&=&u_2(z)\label{2nd}
\end{eqnarray}

The Eq. (\ref{1st}) means that in $D$, $f(z)$ is equal to $0$ if and
only if $z$ is $0$ while the Eq. (\ref{2nd}) implies that $-\delta$
is not in the range of $f(z)$ in $D$. However, $f(z)$ satisfies the
hypothesis of $\displaystyle\frac{1}{16}-$ Theorem, but
$0<|\delta|<\displaystyle\frac{1}{16}$. This is impossible, a
contradiction is obtained and Step 1 is proved.

{\bf Step 2}\quad When $\delta =0$, the three plants defined by
Eq.(\ref{dp}) degenerate to a plant
$p_{1d}(z)=p_{2d}(z)=p_{3d}(z)=0$. Apparently, $c(z)$ meets the
requirement if $c(z)$ is stable.

When $|\delta|=\displaystyle\frac {1}{16}$, since the roots of a
polynomial continuously  depend on its coefficients
\cite{G59,YZH96}, it is easy to know that there does not exist the
desired controller, namely, the three plants defined by Eq.
(\ref{dp}) are not simultaneously stabilzable in this case.

{\bf Step 3}\quad We prove that the three plants defined by
Eq.(\ref{dp}) are simultaneously stabilizable in discrete time case
when $|\delta|
>\displaystyle\frac{1}{16}$.

If $|\delta| >\displaystyle\frac{1}{16}$, from Lemma 1, there exists
an analytic function $f(z)$ on $D$ such that $f(z)=0 \Leftrightarrow
z=0$ and $f'(0)=1$ but $\displaystyle\frac{1}{16} \notin f(D)$.
Moreover, $f(z)$ can be a real analytic function, i.e.,
$f(\overline{z})=\overline{f(z)}$. \footnote{Special constructions
of such functions are presented in \cite{G69,N52}, for instance,
$$f(z)=z\prod_{n=1}^\infty (\displaystyle\frac{1+z^{2n}}{1+z^{2n-1}})^8.$$}

Set
\begin{equation}
g(z)=-16\delta f(\displaystyle\frac{-z}{16\delta}).
\end{equation}

According to the properties of $f(z)$, it is obvious that $g(z)$
satisfies the following properties:

1) $g(\overline{z})=\overline{g(z)}$,

2) $g({z})$ ia analytic on $|z|<16|\delta|$,

3) $g(z)=0\Leftrightarrow z=0$, and $g'(0)=1$,

4) $g({z})$ leaves out $-\delta$ on $|z|<16|\delta|$.

Let
\begin{equation}
\mu = \inf_{z\in \mbox{cl}(D)} |g(z)+\delta|.
\end{equation}
and
\begin{equation}
h(z)=\displaystyle\frac{g(z)-z}{2\delta z^2}.
\end{equation}

From the above properties of $g(z)$, we have $\mu
>0$, $h(z)$ is analytic on $|z|<16|\delta|$, and
\begin{equation}
1+2\delta z h(z)\neq 0, \quad z\in \mbox{cl}(D). \label{h}
\end{equation}

By Runge's Theorem, there exists a real polynomial $q(z)$ such
that
\begin{equation}
|h(z)- q(z)| < \displaystyle\frac{\mu}{2|\delta|}
(\displaystyle\frac{1}{16\delta})^2, \quad z \in \mbox{cl}(D).
\end{equation}

Since $h(z)$ satisfies (\ref{h}), $q(z)$ also satisfies
\begin{equation}
1+2\delta z q(z)\neq 0, \quad z\in \mbox{cl}(D),
\end{equation}
i.e.,
\begin{equation}
1+2\delta z q(z)=u_1(z)\in U(\mbox{cl}D). \label{u1}
\end{equation}

By the definition of $h(z)$, we have
\begin{equation}
|g(z)- z - 2\delta z^2 q(z)| < \mu, \quad z \in \mbox{cl}(D).
\end{equation}

Set
\begin{equation}
p(z)= z + 2\delta z^2 q(z).
\end{equation}
Obviously, $p(0)=0,p'(0)=1$. By the definition of $\mu$,
\begin{equation}
p(z) \neq -\delta, \quad z \in \mbox{cl}(D).
\end{equation}
Hence, we have
\begin{equation}
\delta+z+2\delta z^2 q(z)=u_2(z)\in U(\mbox{cl}D). \label{u2}
\end{equation}

Let $c(z)=q(z)$, by (\ref{u1}) and (\ref{u2}), $c(z)$ is the desired
controller. That completes the proof of Step 3.

Theorem 1 is now proved by combining Steps 1-3.

{\bf Remark 1} \quad {\it Theorem 1 contains a complete
theoretical answer to the `Generalized Champagne Problem "
\cite{P99,PDV02}, automatically including the solution to the
``Champagne Problem" \cite{BG93,BGMR94,BSVW99,P99}. The main
result in Patel \cite{P99} is also contained in Step 1 of the
proof of Theorem 1. Apparently, the proof in this paper is much
more concise and straightforward than that in Patel \cite{P99}.}

{\bf  Remark 2} \quad {\it As a byproduct, Step 2 of the proof of
Theorem 1 contains an answer to an open problem proposed in Patel
\cite{P99}.}

By a similar argument as in the proof of Theorem 1, we have

{\bf Theorem 2} \quad {\it The necessary and sufficient condition
for the existence of a controller that simultaneously stabilizes
the following three plants in continuous time case:
\begin{equation}
p_1(s)=\displaystyle\frac{s-1}{s+1}, \quad
p_2(s)=\displaystyle\frac{(s-1)^2}{((1-\delta )s^2-2\delta
s-(1+\delta)}, \quad p_3(s)=0,\label{bp}
\end{equation}
where $\delta$ is a real number, or, equivalently, simultaneously
stabilizes the following three plants in discrete time case:
\begin{equation}\label{dp2}
p_{1d}(z)= z,\quad p_{2d}(z)=\displaystyle\frac{z^2}{z-\delta
},\quad p_{3d}(z)=0,
\end{equation}
where $\delta$ is a real number, is $\delta =0$ or $|\delta| >
\displaystyle\frac {1}{16}.$}

{\bf Remark 3} \quad {\it Some partial results of Theorem 2 are
available in \cite{B94,BGMR92,BGMR94}. In \cite{B94}, a sufficient
condition was given on the nonexistence of a simultaneously
stabilizing controller for the three plants (in continuous time
case) in Theorem 2.}

{\bf Remark 4} \quad {\it It should be pointed out that Theorem 2 is
equivalent to Theorem 1 in some sense. In fact, a stabilizing
controller multiplied by a non-zero constant is also a stabilizing
controller. When $\delta \neq 0$, by multiplying the three plants in
(\ref{dp2}) by $-2\delta$, and letting $\delta_1=-\delta$, it is
obvious that the three plants obtained are the same as the three
plants in Theorem 1.}

{\bf Remark 5} \quad {\it Although a complete theoretical solution
to the ``Generalized Champagne Problem" is provided in Theorem 1, it
can be seen from the proof of Theorem 1 that it is still difficult
in practice to construct the simultaneously stabilizing controller.
This difficulty can also be seen by the controller design examples
in \cite{LKZ99,PDV02}. A novel stabilizing controller design method
is established in the next section, and the numerical examples
obtained significantly improve the relevant results in
\cite{LKZ99,PDV02}.}

\vskip 0.7cm \section{Controller Design}

Although the ``Generalized Champagne Problem" has been theoretically
resolved, it is still difficult in practice to construct the
simultaneously stabilizing controller. This fact can also be seen in
the specific controller design examples in \cite{LKZ99,PDV02}. In
particular, how can one construct the simultaneously stabilizing
controller with minimal degree for $\delta >
\displaystyle\frac{1}{16}$? And when the degree of the controller is
fixed, how to determine the range of $\delta$? These issues are of
both theoretical and practical significance.

When the degree of the controller is fixed, the controller design
problem of simultaneous stabilization can be transformed in essence
to the problem of how to solve a set of algebraic inequalities.
Early in 1950s Tarski published the well-known work \cite{T51} on
the decidability \cite{ABJ75,T51,W77,W84,W03,YZH96} of this kind of
problems. Tarski's decision algorithm is of theoretical significance
only, since it can not be used to verify any non-trivial algebraic
or geometric propositions in practice \cite{W77,W84,W03,YZH96} due
to its very high computational complexity. The Cylindrical Algebraic
Decomposition (CAD) \cite{ACM184,ACM284,CH91} algorithm proposed by
Collins and subsequently improved by him and his collaborators is
the first practical decision algorithm and can be used to verify
non-trivial algebraic or geometric propositions on computer.
Although, as a generic program for automated theorem proving, its
computational complexity was still very high, the CAD and its
improved variations have become one of the major tools for solving
this kind of problems.

Wu \cite{W77} proposed in 1977 a new decision procedure for proving
geometry theorems. As an important progress in automated theorem
proving \cite{W84,W03,YZH96}, Wu's method is very efficient for
mechanically proving elementary geometry theorems of equality type.
The success of Wu's method has inspired in the world the research of
algebraic approach to automated theorem proving
\cite{W84,W03,YZH96}. However, automated inequality proving has been
a difficult topic in the area of automated reasoning for many years
since the relevant algorithms depend on real algebra and real
geometry. In 1996, Yang and his colleagues \cite{YHZ96,YZH96}
introduced a powerful tool, the Complete Discrimination Systems
(CDS) of Polynomials, for automated reasoning in real algebra. By
means of CDS, many inequality-type theorems from various
applications have been proved or discovered. 

In recent years, by combining discriminant sequences of polynomials
\cite{YHZ96,YZH96} with Wu's method \cite{W77,W84,W03} as well as
Partial Cylindrical Algebraic Decomposition
\cite{ACM184,ACM284,CH91}, Yang et al. proposed some algorithms
which are able to discover new inequalities. These algorithms are
complete for an extensive class of problems involving inequalities
and are applicable to the controller design in simultaneous
stabilization. Based on these algorithms, a generic program called
``Discoverer" \cite{YHZ96,YX05} and a generic program called
``Bottema" \cite{Y98,Y99,YX03} were implemented as Maple packages.
The following controller examples are all obtained automatically by
calling ``Discoverer" \cite{YHZ96,YX05} or ``Bottema"
\cite{Y98,Y99,YX03}.

The following sufficient condition for checking the Hurwitz
stability of a polynomial may be used for improving computational
efficiency when calling the above programs.

{\bf Lemma 3} \cite{N76} \quad {\it Suppose $f(s)=a_0+a_1s+\cdots
+a_ns^n~ (a_i>0,i=0,1,2,\cdots,n,n\geq 3)$, $f(s)$ is Hurwitz stable
if
\begin{equation}
a_{i-1}a_{i+2}\leq 0.4655a_ia_{i+1},\quad (i=1,2,\cdots,n-2).
\end{equation}}

In addition, to get a proper controller, we introduce the following
theorem on the disturbance to the Hurwitz polynomials.

{\bf Theorem 3} \quad {\it Given a polynomial $f(s)=a_0+a_1s+\cdots
+a_ms^m~ (a_i>0,i=0,1,2,\cdots,m) \in H^m$, then for any integer
$n>m$ and sufficiently small $\varepsilon>0$, there exist
$\varepsilon_i$, $0<\varepsilon_i<\varepsilon
~(i=m+1,m+2,\cdots,n)$, such that $g(s)=a_0+a_1s+\cdots +a_ms^m
+\varepsilon_{m+1}s^{m+1}+\varepsilon_{m+2}s^{m+2}+\cdots+\varepsilon_ns^n\in
H^n.$}

{\bf Proof} \quad Suppose $f(s)=a_0+a_1s+\cdots +a_ms^m
(a_i>0,i=0,1,2,\cdots,m) \in H^m$, we only need to prove that
$g(s)=a_0+a_1s+\cdots +a_ms^m+\varepsilon s^{m+1}\in H^{m+1}$ for
sufficiently small $\varepsilon>0$ because we can treat $g(s)$ with
the same manner iteratively.

Since $f(s)=a_0+a_1s+\cdots +a_ms^m\in H^m$, it is well known that
Hurwitz principal minors $\Delta^f_i>0~ (i=1,2,\cdots,m)$
\cite{G59}.

By computing Hurwitz principal minors of $g(s)=a_0+a_1s+\cdots
+a_ms^m+\varepsilon s^{m+1}$, it is easy to verify that
$\Delta^g_1=a_m, \Delta^g_i = a_m\Delta^f_{i-1}
 + \varepsilon\delta_i(\varepsilon,a_0,a_1,\cdots,a_m)
 ~(i=2,\cdots,m,m+1)$,
where $\delta_i(\varepsilon,a_0,a_1,\cdots,a_m) ~
(i=2,\cdots,m,m+1)$ are polynomials in
$(\varepsilon,a_0,a_1,\cdots,a_m)$. Obviously, when $\varepsilon>0$
is sufficiently small, we have $\Delta^g_i>0 ~
(i=1,2,\cdots,m,m+1)$. This completes the proof.

In what follows, we will provide some controller examples for
``Generalized Champagne Problem" of simultaneous stabilization by
applying ``Discoverer" \cite{YHZ96,YX05} or ``Bottema"
\cite{Y98,Y99,YX03}. Without loss of generality, we suppose the
denominator polynomials of the controllers are monic. In addition,
for convenience, we only consider the case $\delta>0$.

{\bf Example 1} \quad {\it For the three plants defined in
``Generalized Champagne Problem", when the degree of controllers is
restricted to $0$, i.e., the controllers are constants, invoking
``Discoverer" or ``Bottema", we obtain the following results. There
exist desired controllers for $\delta>\displaystyle\frac{1}{2}$,
whereas there do not exist requested controllers for $\delta\leq
\displaystyle\frac{1}{2}$. For a given
$\delta>\displaystyle\frac{1}{2}$, for example,
$\delta=\displaystyle\frac{3}{4}$, $c(s)=y_0$ is a desired
controller if and only if
$\displaystyle\frac{1}{6}<y_0<\displaystyle\frac{1}{2}$.}

{\bf Remark 6} \quad {\it Example 1 looks simple, but it shows
that: 1) When the degree of the controller is restricted, the
explicit bound of $\delta$ can be obtained by using the packages
developed by Yang et al.; 2) Moreover, when $\delta$ is fixed, the
ranges of the parameters of controllers can be obtained. Due to
the completeness of these algorithms, the conditions obtained are
both necessary and sufficient.}

{\bf Example 2} \quad {\it For the three plants defined in
``Generalized Champagne Problem", when the numerator polynomials of
controllers are restricted to $0-th$ order polynomials and the
denominator polynomials of controllers to $1st$ order polynomials,
i.e., the controllers have the form
$c(s)=\displaystyle\frac{y_0}{s+x_0}$, invoking ``Discoverer" or
``Bottema", we have the following results as in the above example.
There exist desired controllers for
$\delta>\displaystyle\frac{1}{2}$, whereas there do not exist
requested controllers for $\delta\leq \displaystyle\frac{1}{2}$. For
a given $\delta>\displaystyle\frac{1}{2}$, say,
$\delta=\displaystyle\frac{3}{4}$, $c(s)=\displaystyle\frac{1}{s+3}$
is a desired controller.}

{\bf Example 3} \quad {\it  For the three plants defined in
``Generalized Champagne Problem", when the numerator polynomials of
controllers are restricted to $1st$ order polynomials and the
denominator polynomials of controllers to $0-th$ order polynomials,
i.e., the controllers have the form  $c(s)=y_1s+y_0$, invoking
``Discoverer" or ``Bottema", we have the following results. There
exist desired controllers for $\delta>\displaystyle\frac{1}{4}$,
whereas there do not exist requested controllers for $\delta\leq
\displaystyle\frac{1}{4}$. For a given
$\delta>\displaystyle\frac{1}{4}$, say
$\delta=\displaystyle\frac{1}{2}$, if $(y_1,y_0)\in
\{(\displaystyle\frac{1}{10},\displaystyle\frac{51}{100}),
(\displaystyle\frac{3}{5},\displaystyle\frac{3}{5}),
(\displaystyle\frac{4}{5},\displaystyle\frac{3}{5})\}$,
$c(s)=y_1s+y_0$ is a requested controller. Moreover, to get a proper
controller, by the continuous dependance for roots of polynomials on
their coefficients \cite{G59,YZH96} and Theorem 3, if
$(x_1,y_1,y_0)\in
\{(\varepsilon,\displaystyle\frac{1}{10},\displaystyle\frac{51}{100}),
(\varepsilon,\displaystyle\frac{3}{5},\displaystyle\frac{3}{5}),
(\varepsilon,\displaystyle\frac{4}{5},\displaystyle\frac{3}{5})\}$
and $\varepsilon>0$ is sufficiently small, e.g.,
$\varepsilon=\displaystyle\frac{1}{10}$,
$\widetilde{c}(s)=\displaystyle\frac{y_1s+y_0}{\varepsilon s+1}$ is
a desired proper controller.}

{\bf Example 4} \quad {\it For the three plants defined in
``Generalized Champagne Problem", when the numerator polynomials of
controllers are restricted to $1st$ order polynomials and the
denominator polynomials of controllers to $1st$ order polynomials,
i.e., the controllers have the form
$c(s)=\displaystyle\frac{y_1s+y_0}{s+x_0}$, invoking ``Discoverer"
or ``Bottema", we have that there exist desired controllers for
$\delta>\displaystyle\frac{1}{4}$ whereas there do not exist
requested controllers for $\delta\leq \displaystyle\frac{1}{4}$. For
a given $\delta>\displaystyle\frac{1}{4}$, say,
$\delta=\displaystyle\frac{1}{3}$, if $(x_0,y_1,y_0)\in \{
(2,{\displaystyle\frac {31}{10}},{\displaystyle\frac
{201}{100}}),(3,\displaystyle\frac{41}{10},{\displaystyle\frac
{301}{100}}), (6,8,{\displaystyle\frac {61}{10}})\}$,
$c(s)=\displaystyle\frac{y_1s+y_0}{s+x_0}$ is a desired controller.}

{\bf Remark 7} \quad {\it The controllers obtained in the above
examples are the sample points picked out from the cells of some
appropriate decomposition of the parametric space which satisfy the
requirements of simultaneous stabilization. As mentioned in Example
1, when the degrees of controllers are restricted, the explicit
bound of $\delta$ and the ranges of the parameters of the
controllers can be obtained by further analysis on the distribution
of the roots of the polynomials output by ``Discoverer".}

{\bf Example 5} \quad {\it For the three plants defined in
``Generalized Champagne Problem", when the numerator polynomials of
controllers are restricted to $2nd$ order polynomials and the
denominator polynomials of controllers to $0-th$ order polynomials,
i.e., the controllers have the form $c(s)=y_2s^2+y_1s+y_0$, invoking
``Discoverer", we have that there do not exist desired controllers
for $\delta=\displaystyle\frac{1}{6}$. When
$\delta=\displaystyle\frac{10}{59}$, 
using the sufficient condition proposed in Lemma 3 for Hurwitz
stability, we have that there exist such kind of controllers. For
example, if $(y_2,y_1,y_0)\in \{ $ $({\displaystyle\frac
{191}{100}}, $ ${\displaystyle\frac {39001}{10000}}, $
${\displaystyle\frac {2450001}{1000000}}),$ $({\displaystyle\frac
{191}{100}}, $ ${\displaystyle\frac {390019}{100000}}, $
${\displaystyle\frac {2450003}{1000000}}),$ $({\displaystyle\frac
{97}{50}}, $ ${\displaystyle\frac {39001}{10000}}, $
${\displaystyle\frac {245001}{100000}}),$ $({\displaystyle\frac
{97}{50}}, $ ${\displaystyle\frac {19501}{5000}}, $
${\displaystyle\frac {245001}{100000}}),$ $({\displaystyle\frac
{19501}{10000}}, $ ${\displaystyle\frac {39003}{10000}}, $
${\displaystyle\frac {4900299}{2000000}})\}$, $c(s)=y_2s^2+y_1s+y_0$
is a desired controller. Moreover, to get a proper controller, by
the continuous dependance for roots of polynomials on their
coefficients \cite{G59,YZH96} and Theorem 3, if
$(x_2,x_1,y_2,y_1,y_0)\in \{ (\varepsilon, $ $\varepsilon, $
${\displaystyle\frac {191}{100}}, $ ${\displaystyle\frac
{39001}{10000}}, $ ${\displaystyle\frac {2450001}{1000000}}),$
$(\varepsilon, $ $\varepsilon, $ ${\displaystyle\frac {191}{100}}, $
${\displaystyle\frac {390019}{100000}}, $ ${\displaystyle\frac
{2450003}{1000000}}),$ $(\varepsilon, $ $\varepsilon, $
${\displaystyle\frac {97}{50}}, $ ${\displaystyle\frac
{39001}{10000}}, $ ${\displaystyle\frac {245001}{100000}}),$
$(\varepsilon, $ $\varepsilon, $ ${\displaystyle\frac {97}{50}}, $
${\displaystyle\frac {19501}{5000}}, $ ${\displaystyle\frac
{245001}{100000}}),$ $(\varepsilon, $ $\varepsilon, $
${\displaystyle\frac {19501}{10000}}, $ ${\displaystyle\frac
{39003}{10000}}, $ ${\displaystyle\frac {4900299}{2000000}})\}$ and
$\varepsilon>0$ is sufficiently small, e.g.,
$\varepsilon=\displaystyle\frac{1}{10000000}$, then
$\widetilde{c}(s)=\displaystyle\frac{y_2s^2+y_1s+y_0}{\varepsilon
s^2+\varepsilon s+1}$ is a desired proper controller.}

{\bf Remark 8} \quad {\it The value of $\delta$ obtained in Example
5 is an improvement over the bound conjectured in \cite{LKZ99}.
Although an improvement over this bound was also made by Patel et
al. \cite{PDV02}, the degree of the controller provided by them
equals $9$ whereas the controller presented in Example 5 is of
degree $2$ which is much more lower. In addition, we can see from
the above examples that the improvement over the bound of $\delta$
conjectured in \cite{LKZ99} could not be achieved by controllers
with degrees less than $2$.}

{\bf Example 6} \quad {\it For the three plants defined in
``Generalized Champagne Problem", suppose the numerator polynomials
of controllers are restricted to $3rd$ order polynomials and the
denominator polynomials of controllers to $0-th$ order polynomials,
i.e., the controllers have the form $c(s)=y_3 s^3+y_2s^2+y_1s+y_0$.
When $\delta=\displaystyle\frac{1}{7}$, using the sufficient
condition proposed in Lemma 3 for Hurwitz stability, we have that
there exist such kind of controllers. For example, if
$(y_3,y_2,y_1,y_0)\in \{ $ $({\displaystyle\frac {1037}{1000}}, $
${\displaystyle\frac {30077}{10000}}, $ ${\displaystyle\frac
{50001}{10000}}, $ ${\displaystyle\frac {300001}{100000}}),$
$({\displaystyle\frac {26}{25}}, $ ${\displaystyle\frac {376}{125}},
$ ${\displaystyle\frac {50001}{10000}}, $ ${\displaystyle\frac
{300001}{100000}}),$ $({\displaystyle\frac {113}{100}}, $
${\displaystyle\frac {378}{125}}, $ ${\displaystyle\frac
{50001}{10000}}, $ ${\displaystyle\frac {300001}{100000}}),$
$({\displaystyle\frac {57029}{50000}}, $ ${\displaystyle\frac
{121}{40}}, $ ${\displaystyle\frac {50001}{10000}}, $
${\displaystyle\frac {300001}{100000}}),$ $({\displaystyle\frac
{11407}{10000}}, $ ${\displaystyle\frac {121}{40}}, $
${\displaystyle\frac {50001}{10000}}, $ ${\displaystyle\frac
{300001}{100000}})\},$ $({\displaystyle\frac {114113}{100000}}, $
${\displaystyle\frac {1513}{500}}, $ ${\displaystyle\frac
{50001}{10000}}, $ ${\displaystyle\frac {300001}{100000}})\}$,
$c(s)=y_3 s^3+y_2s^2+y_1s+y_0$ is a desired controller. Moreover, to
get a proper controller, by the continuous dependance for roots of
polynomials on their coefficients \cite{G59,YZH96} and Theorem 3, if
$(x_3,x_2,x_1,y_3,y_2,y_1,y_0)\in $ $\{(\varepsilon_1, $
$\varepsilon, $ $\varepsilon, $ ${\displaystyle\frac {1037}{1000}},
$ ${\displaystyle\frac {30077}{10000}}, $ ${\displaystyle\frac
{50001}{10000}}, $ ${\displaystyle\frac {300001}{100000}}),$
$(\varepsilon_1, $ $\varepsilon,$ $\varepsilon, $
${\displaystyle\frac {26}{25}}, $ ${\displaystyle\frac {376}{125}},
$ ${\displaystyle\frac {50001}{10000}}, $ ${\displaystyle\frac
{300001}{100000}}),$ $(\varepsilon_1, $ $\varepsilon, $
$\varepsilon, $ ${\displaystyle\frac {113}{100}}, $
${\displaystyle\frac {378}{125}}, $ ${\displaystyle\frac
{50001}{10000}}, $ ${\displaystyle\frac {300001}{100000}}),$
$(\varepsilon_1, $ $\varepsilon, $ $\varepsilon, $
${\displaystyle\frac {57029}{50000}}, $ ${\displaystyle\frac
{121}{40}}, $ ${\displaystyle\frac {50001}{10000}}, $
${\displaystyle\frac {300001}{100000}}),$ $(\varepsilon_1, $
$\varepsilon, $ $\varepsilon, $ ${\displaystyle\frac
{11407}{10000}}, $ ${\displaystyle\frac {121}{40}}, $
${\displaystyle\frac {50001}{10000}}, $ ${\displaystyle\frac
{300001}{100000}}),$ $(\varepsilon_1, $ $\varepsilon, \varepsilon, $
${\displaystyle\frac {114113}{100000}}, $ ${\displaystyle\frac
{1513}{500}}, $ ${\displaystyle\frac {50001}{10000}}, $
${\displaystyle\frac {300001}{100000}})\}$ and
$\varepsilon^2>\varepsilon_1$ for sufficiently small
$\varepsilon_1>0$ and $\varepsilon>0$, e.g., $\varepsilon=10^{-7}$
and $\varepsilon_1=10^{-15}$,
$\widetilde{c}(s)=\displaystyle\frac{y_3s^3+y_2s^2+y_1s+y_0}{\varepsilon_1
s^3+\varepsilon s^2+\varepsilon s+1}$ is a desired proper
controller.}

{\bf Remark 9} \quad {\it The value of $\delta$ appeared in Example
6 is an improvement over the minimal bound proposed in \cite{PDV02}.
The degree of the controller presented in \cite{PDV02} equals 9
whereas the controller given in Example 6 is of degree $3$ which is
much more lower.}

{\bf Remark 10} \quad {\it The above algorithms are also complete
for high-order controller design. From our computational
experiments, other phenomena can be observed. For instance, the
improvement on the bound of $\delta$ in ``Generalized Champagne
Problem" mainly depends on the increase on the order of numerator
polynomial of the stabilizing controller. These problems deserve
further research and are omitted here for succinctness.}

{\bf Remark 11} \quad {\it The generic programs ``Discoverer"
\cite{YHX01,YX05} and ``Bottema" \cite{Y98,Y99,YX03} are powerful
tools in practice for automated inequalities proving and can be
applied in various fields. In this paper, only some basic functions
of those packages are employed to determine the ranges of parameters
and to design simultaneously stabilizing controllers. In previous
examples, it is demonstrated how powerful these packages are.
Indeed, according to the {\em Finiteness Theorem} proposed by W.T.
Wu \cite{W94} on global-optimization, global-optimization problems
can be theoretically solved by `Discoverer" \cite{YHX01,YX05} and
``Bottema" \cite{Y98,Y99,YX03} where the objectives are polynomials
and constraints are also polynomial equalities or inequalities. The
potential applications of the packages are considerable.}

{\bf Remark 12} \quad {\it  It should be pointed out that, although
``Discoverer" \cite{YHX01,YX05} and Bottema \cite{Y98,Y99,YX03} are
powerful, the computational complexity increases very quickly with
the dimension, i.e., the number of parameters. It is a problem
deserving further study how to promote the efficiencies of these
algorithms for dealing with the large-scale engineering optimization
problems. Combining symbolic computation with numerical calculation
as well as large-scale parallel numerical processing may be an
effective way \cite{W03}.}

\section{Conclusion}

In this paper, the well-known ``Generalized Champagne Problem" for
simultaneous stabilization of linear systems has been resolved by
using complex analysis \cite{A73,C78,G69,N52,R87} and Blondel's
technique \cite{B94,BG93,BGMR94}. We gave a complete answer to the
open problem proposed in Patel et al. \cite{P99,PDV02}, which
automatically includes the solution to the original ``Champagne
Problem" \cite{BG93,BGMR94, BSVW99,LKZ99,P99,PDV02}. Based on the
recent developments in automated inequality-type theorem proving
\cite{Y98,Y99,YHX01,YX03,YX05}, a novel stabilizing controller
design method has been established. Our numerical examples
significantly improved the relevant results in the literature
\cite{LKZ99,PDV02}.

\section*{Acknowledgements}

The authors express their thanks to Prof Wenjun Yuan of Guangzhou
Univ and Dr Nong Gu of Dekin Univ in Australia for providing some
relevant papers.

\vskip 20pt
\vspace*{1\baselineskip}

\end{document}